\documentclass[12pt,reqno]{amsart}

\usepackage{hyperref}

\usepackage{amsthm,amssymb,latexsym}
\usepackage{psfig,epsf}
\usepackage{float,color}

\newcommand{\bega}{\begin{eqnarray}}
\newcommand{\ega}{\end{eqnarray}}
\newcommand{\bb}{\begin{equation}}
\newcommand{\ee}{\end{equation}}

 \newtheorem{thm}{Theorem}[section]
 
 \newtheorem{lem}{Lemma}[section]
 \newtheorem{cor}{Corollary}[section]
 \newtheorem{exm}{Example}[section]

\allowdisplaybreaks[1]
\begin{document}

\begin{center}
\epsfxsize=4in
%\leavevmode\epsffile{logo130.eps}
\end{center}

\medskip

\begin{center}
{\Large{ The Hankel Transform of the Sum of \\
\medskip Consecutive Generalized Catalan Numbers }}
\end{center}

\medskip

\begin{center}

{\bf  Predrag Rajkovi\'c, Marko D. Petkovi\'c,}\\
 University of Ni\v s, Serbia and Montenegro  \\

\smallskip

{\bf Paul Barry} \\
School of Science, Waterford Institute of Technology, Ireland
\end{center}

%\maketitle

\medskip

%\begin{center}
{\bf Abstract.} {\em We discuss the properties of the Hankel
transformation of a sequence whose elements are the sums of
consecutive generalized Catalan numbers and find their values
in the closed form. }
%\end{center}

\medskip
% 32A05 Power series, series of functions
% 11Y55 Calculation of integer sequences

% 33C47 Other special orthogonal polynomials and functions
% 34A25 Analytical theory: series, transformations, transforms,
% operational calculus, etc. [See also 44-xx]

\smallskip


 \section{Introduction}

The {\it Hankel transform} of a given sequence
$A=\{a_0,a_1,a_2,...\}$ is the sequence of Hankel determinants
$\{h_0, h_1, h_2,\dots \}$ (see Layman \cite{layman}) where
$h_{n}=|a_{i+j-2}|_{i,j=1}^{n}$, i.e

%\begin{center}
\bb
 \label{gen1}
 A=\{a_n\}_{n\in\mathbb N_0}\quad \rightarrow \quad
 h=\{h_n\}_{n\in\mathbb N_0}:\quad
h_n=\left|
\begin{array}{ccccc}
 a_0\ & a_1\  & \cdots & a_n  &  \\
 a_1\ & a_2\  &        & a_{n+1}  \\
\vdots &      & \ddots &          \\
 a_n\ & a_{n+1}\ &    & a_{2n}
\end{array}
\right|
\ee
%\end{center}

In this paper, we will consider the sequence of  the sums of two
adjacent generalized Catalan numbers with parameter $L$:
%\begin{eqnarray*}
\bb
 \label{gen2}
 a_0=L+1, \qquad  a_n=a_n(L)=c(n;L)+c(n+1;L)  \qquad (n\in \mathbb N),
\ee
where
\bb
 \label{gen3}
 c(n;L)=T(2n,n;L)-T(2n,n-1;L),
\ee
% \end{eqnarray*}
with
%\begin{eqnarray*}
\bb
 \label{gen}
 T(n,k;L)=\sum_{j=0}^{n-k} \binom{k}{j} \binom{n-k}{j}L^j.
\ee
%\end{eqnarray*}

%{\bf Example 1.1.}

\begin{exm}\rm
Let $L=1$.
Vandermonde's convolution identity implies that
$$
\binom{n}{k} = \sum_{j} \binom{k}{j} \binom{n-k}{j}.
$$
Hence
$$
T(2n,n;1)= \binom{2n}{n}, \qquad T(2n,n-1;1)= \binom{2n}{n-1},
$$
wherefrom we get Catalan numbers
$$
c(n)=\binom{2n}{n} - \binom{2n}{n-1} = \frac1{n+1} \binom{2n}{n}
$$
and
\begin{eqnarray*}
  a_n=c(n)+c(n+1)=
\frac{(2n)!(5n+4)}{n!(n+2)!}  \qquad (n=0,1,2,\ldots).
\end{eqnarray*}

In the paper \cite{CRI}, A. Cvetkovi\'c, P. Rajkovi\'c and M.
Ivkovi\'c have proved that the Hankel transform of $a_n$ equals
sequence of Fibonacci numbers with odd indices
$$
h_n=F_{2n+1}=\frac{1}{ \sqrt5 \
2^{n+1}}\left\{(\sqrt5+1)(3+\sqrt5)^n+(\sqrt5-1)(3-\sqrt5)^n\right\}.
$$

\end{exm}

\smallskip

%{\bf Example 1.2.}
\begin{exm}\rm
\ For $L=2$ we get like $a_n(2)$ the next
numbers
$$
3,8,28,112,484,\ldots ,
$$
and the Hankel transform $h_n$:
$$
3,20,272,7424,405504,\ldots .
$$
One of us,  P. Barry conjectured that
$$
h_n(2)=2^{\frac{n^2-n}2-2}\left\{(2+\sqrt2)^{n+1}+(2-\sqrt2)^{n+1}\right\}.
$$
\end{exm}

In general,  P. Barry made the conjecture, which we will prove
through this paper.

\begin{thm} {\bf (The main result)} For the generalized Pascal triangle associated to the sequence $n\mapsto L^n$,
the Hankel transform of the sequence
$$
c(n;L)+c(n+1;L)
$$
is given by
\begin{eqnarray}
\aligned &h_n= \frac{L^{(n^2-n)/2}}{2^{n+1}\sqrt{L^2+4}}\cdot\\
& \left\{(\sqrt{L^2+4}+L)(\sqrt{L^2+4}+L+2)^{n}
+(\sqrt{L^2+4}-L)(L+2-\sqrt{L^2+4})^{n}\right\}.
\endaligned
\end{eqnarray}
\end{thm}

From now till the end, let us denote by
\begin{eqnarray}
\xi=\sqrt{L^2+4},\quad t_1=L+2+\xi, \quad t_2=L+2-\xi.
\end{eqnarray}
Now, we can write
$$
h_n= \frac{L^{n(n-1)/2}}{2^{n+1}\xi}\cdot \left((\xi+L)t_1^{n}
+(\xi-L)t_2^{n}\right).
$$
Or, introducing
\begin{eqnarray}
\varphi_n = t_1^{n}+t_2^{n}, \quad \psi_n = t_1^{n}-t_2^{n} \qquad
(n\in\mathbb N_0),
\end{eqnarray}
the final statement can be expressed by
\begin{eqnarray}
h_n= \frac{L^{n(n-1)/2}}{2^{n+1}\xi}\cdot \left(L\psi_n+\xi
\varphi_n\right).
\end{eqnarray}

\begin{lem} The values $\varphi_n$ and $\psi_n$ satisfy the next
relations
\begin{eqnarray}
\varphi_j\cdot \varphi_k= \varphi_{j+k} +  (4L)^j\varphi_{k-j},
\quad \psi_j\cdot \psi_k= \varphi_{j+k} -
(4L)^j\varphi_{k-j}\qquad (0\leq j \leq k)
\end{eqnarray}
\begin{eqnarray}
\varphi_j\cdot \psi_k= \psi_{j+k} +  (4L)^j\psi_{k-j}, \qquad
\psi_j\cdot \varphi_k= \psi_{j+k} - (4L)^j\psi_{k-j}\qquad (0\leq
j \leq k).
\end{eqnarray}
\end{lem}

 \begin{cor}
Assuming that the main theorem is true, the function $h_n=h_n(L)$ is the next polynomial
$$
\aligned
h_n(L) &= 2^{-n} L^{n(n-1)/2}\\
&  \cdot\left\{\sum_{i=0}^{[(n-1)/2]}
\tbinom{n}{2i+1}L(L+2)^{n-2i-1}(L^2+4)^i + \sum_{i=0}^{[n/2]}
\tbinom{n}{2i}(L+2)^{n-2i}(L^2+4)^i\right\}.
\endaligned
$$
\end{cor}
\noindent {\sl Proof.} By previous notation, we can write
$$
\aligned
&(L+\xi)(L+2+\xi)^n - (L-\xi)(L+2-\xi)^n\\
&=(L+\xi)\sum_{k=0}^n \tbinom{n}{k}(L+2)^{n-k}\xi^k
- (L-\xi)\sum_{k=0}^n (-1)^k \tbinom{n}{k}(L+2)^{n-k}\xi^k\\
&= \sum_{k=0}^n \bigl(1-(-1)^k\bigr)\tbinom{n}{k}L(L+2)^{n-k}\xi^k
+ \sum_{k=0}^n \bigl(1+(-1)^k\bigr)\tbinom{n}{k}(L+2)^{n-k}\xi^{k+1}  \\
&= 2\sum_{i=0}^{[(n-1)/2]}
\tbinom{n}{2i+1}L(L+2)^{n-2i-1}\xi^{2i+1} + 2\sum_{i=0}^{[n/2]}
\tbinom{n}{2i}(L+2)^{n-2i}\xi^{2i+1} \\
&= 2\xi \left\{\sum_{i=0}^{[(n-1)/2]}
\tbinom{n}{2i+1}L(L+2)^{n-2i-1}\xi^{2i} + \sum_{i=0}^{[n/2]}
\tbinom{n}{2i}(L+2)^{n-2i}\xi^{2i}\right\},
\endaligned
 $$
wherefrom immediately follows the polynomial expression for
$h_n$.$\Box$

% section 2
\section{The generating function for the sequences \\of numbers and orthogonal polynomials}

The Jacobi polynomials are given by
$$
P^{(a,b)}_n(x) = \frac{1}{2^n} \sum_{k=0}^n \binom{n+a}{k}
\binom{n+b}{n-k} (x-1)^{n-k}(x+1)^k\quad (a,b>-1).
$$
Also, they can be written in the form
$$
P^{(a,b)}_n(x) =\Bigl(\frac{x-1}{2}\Bigr)^n \sum_{k=0}^n
\binom{n+a}{k} \binom{n+b}{n-k} \Bigl(\tfrac{x+1}{x-1}\Bigr)^k.
$$
From the fact
$$
L=\frac{x+1}{x-1} \quad \Leftrightarrow \quad
x=\frac{L+1}{L-1}\qquad (x\ne1,\ L\ne1).
$$
we conclude that
$$\aligned
T(2n,n;L)&=(L-1)^n  \cdot P^{(0,0)}_n\Bigl(\tfrac{L+1}{L-1}\Bigr),
\\
 T(2n+2,n;L)&=(L-1)^n  \cdot
P^{(2,0)}_n\Bigl(\tfrac{L+1}{L-1}\Bigr). \endaligned
$$
The  generating function $G(x,t)$ for the Jacobi polynomials is \bb
 \label{gen5}
G^{(a,b)}(x,t) = \sum_{n=0}^{\infty}  P^{(a,b)}_n(x)  t^n =
\frac{2^{a+b}}{\phi \cdot (1-t+\phi)^a \cdot (1+t+\phi)^b},
\ee
where
$$
\phi=\phi(x,t)=\sqrt{1-2xt+t^2}.
$$
Now,
$$
\aligned
 \sum_{n=0}^{\infty}  T(2n,n;L)\  t^n &=\sum_{n=0}^{\infty}
P^{(0,0)}_n\Bigl(\tfrac{L+1}{L-1}\Bigr)  \bigl((L-1)t\bigr)^n =
G^{(0,0)}\Bigl(\tfrac{L+1}{L-1},(L-1)t  \Bigr),\\
\sum_{n=0}^{\infty}  T(2n+2,n;L)\  t^n &=\sum_{n=0}^{\infty}
P^{(2,0)}_n\Bigl(\tfrac{L+1}{L-1}\Bigr)  \bigl((L-1)t\bigr)^n =
G^{(2,0)}\Bigl(\tfrac{L+1}{L-1},(L-1)t  \Bigr).
\endaligned
$$
Also,
$$
\aligned \sum_{n=0}^{\infty}  T(2n,n-1;L)\  t^n  &=t\cdot \left\{
G^{(2,0)}\Bigl(\tfrac{L+1}{L-1},(L-1)t  \Bigr) - 1\right\},\\
\sum_{n=0}^{\infty}  T(2n+2,n+1;L)\  t^n  &=\frac1{t}\cdot \left\{
G^{(0,0)}\Bigl(\tfrac{L+1}{L-1},(L-1)t  \Bigr) - 1\right\}.
\endaligned
$$
The generating function $\mathcal{G}(t;L)$ for the sequence
$\{a_n\}_{n \geq 0}$ is given by
\bb
 \label{gen6}
\aligned
\mathcal{G}(t;L) &= \sum_{n=0}^{\infty}a_n t^n \\
&= \frac{t+1}{t} G^{(0,0)}\Bigl(\tfrac{L+1}{L-1},(L-1)t  \Bigr) -
(t+1) G^{(2,0)}\Bigl(\tfrac{L+1}{L-1},(L-1)t  \Bigr)  - \frac1{t}.
\endaligned
\ee

After some computation, we prove the next theorem.

\begin{thm}
The generating function $\mathcal{G}(t;L)$ for the sequence
$\{a_n\}_{n \geq 0}$ is
\bb
 \label{gen8}
\mathcal{G}(t;L)= \frac{t+1}{\rho(t;L)}\ \left\{\frac{1}{t}
 - \frac{4}{(1-(L-1)t+ \rho(t;L))^2} \right\}
- \frac1{t},
 \ee
where
 \bb
 \label{gen7}
\rho (t;L) = \phi\Bigl(\tfrac{L+1}{L-1},(L-1)t  \Bigr) =
\sqrt{1-2(L+1)t+(L-1)^2t^2}
 \ee
\end{thm}

The function $\rho (t;L)$ has domain
$$
D_ {\rho}= \Bigl(-\infty, \frac{1-2\sqrt{L}+L}{1-2L+L^2}\Bigr)\cup
\Bigl(\frac{1+2\sqrt{L}+L}{1-2L+L^2},+\infty\Bigr)  \qquad (L\ne
1),
$$
and
$$
D_ {\rho}= (-\infty, 1/4) \qquad (L=1).
$$

\smallskip

%{\bf Example 2.1.}
\begin{exm}\rm
\ For $L=1$, we get
\bb
 \label{gen9}
\mathcal{G}(t;1)=\sum_{n=0}^{\infty}a_n(1)\ t^n=\frac1{t} \left(
\frac{(1-\sqrt{1-4t})(1+t)}{2t}-1 \right).
 \ee
and for $L=2$, we find
\bb
 \label{gen10}
\mathcal{G}(t;2)=\sum_{n=0}^{\infty}a_n(2)\ t^n=
-\frac1{t}+\frac{t+1}{\sqrt{t^2-6t+1}} \left\{ \frac1{t}-
\frac{4}{(1-t+\sqrt{t^2-6t+1})^2} \right\} .
 \ee
 \end{exm}

%\section 3
\section{The weight function corresponding to the functional}

It is known (for example, see Krattenthaler \cite{kratt}) that the
Hankel determinant $h_n$ of  order $n$ of the sequence
$\{a_n\}_{n\geq 0}$ equals
 \bb
 \label{formula}
 h_n=a_0^n \beta_1^{n-1} \beta_2^{n-2} \cdots \beta_{n-2}^2 \beta_{n-1},
 \ee
  where $ \{ \beta_n \}_{n \geq 1}$ is the  sequence given by:

\bb \label{bes} \mathcal{G}(x)=\sum_{n=0}^{\infty}a_n x^n=
\frac{a_0}{\displaystyle 1+\alpha_0x -\frac{\displaystyle
\beta_1x^2}{1+\alpha_1x-\frac{\displaystyle
\beta_2x^2}{\displaystyle 1+\alpha_2x- \cdots}}} \ee

The sequences $ \{ \alpha_n \}_{n \geq 0}$ and $ \{ \beta_n \}_{n
\geq 1}$ are the coefficients in the recurrence relation

\bb Q_{n+1}(x)=(x-\alpha_n)Q_n(x) -\beta_n Q_{n-1}(x), \ee

 \noindent where  $\{Q_n(x)\}_{n \geq 0}$ is the monic polynomial sequence  orthogonal with respect to
  the functional $\mathcal U$ determined by

\bb \label{funk} \mathcal U[x^n]=a_n \quad (n=0,1,2,\ldots).
 \ee

In this section the functional will be  constructed for the sum of
consecutive generalized Catalan numbers.

We would like  to express $\mathcal U[f]$ in the form:
$$
\mathcal U[f(x)]=\int_R f(x) d\psi (x),
$$
\noindent where $\psi (x)$ is a distribution, or, even more, to
find the weight function $w(x)$ such that $w(x)=\psi'(x).$

 Denote by $F(z;L)$ the function
  $$
  F(z;L)=\sum_{k=0}^{\infty} a_k z^{-k-1}.
  $$
 From the generating function (\ref{gen8}), we have:
 \bb
 F(z;L)=z^{-1}\ \mathcal{G}\bigl(z^{-1};L ).
 \ee
 and after some simplifications we obtain that
$$
\aligned
 F(z;L)&=-1+\frac{2(z+1)}{L-1+z+\sqrt{L^2+(z-1)^2-2L(z+1)}}\\
       & =-1+\frac{2(z+1)}{L-1+z(1+z\rho(\frac1z,L))}
 \endaligned
 $$

%{\bf Example 3.1.}
 \begin{exm}\rm
 From  (\ref{gen9}) and (\ref{gen10}), we yield
$$
\aligned
 F(z;1)&=z^{-1}\ \mathcal{G}\bigl(z^{-1};1 )=\frac1{2} \left\{ z-1-(z+1)\sqrt{1-\frac4{z}} \right\},\\
 F(z;2)&=\frac{-1}{2z} \left\{  1+z\Biggl(2-z+(z+1)\sqrt{1-\frac6{z}+\frac1{z^2}}\Biggr)\right\}.
 \endaligned
$$
Notice that
$$
\aligned
\int  F(z;2) dz = z &+ \frac14 z(z-1)
\rho(1/z,2)+\log(z)\\
&-\frac12 \log\Bigl(1+z\bigl(\rho(1/z,2)-3\bigr)\Bigr)-\frac72
\log\bigl(z-3+z\rho(1/z,2)\bigr).
\endaligned
$$
It will be the impulse for further discussion.
 \end{exm}

 Denote by
 $$
 R(z;L)=z\rho(\frac1z,L)=\sqrt{L^2+(z-1)^2-2L(z+1)}.
 $$
 From the theory of distribution functions (see Chihara
  \cite{chihara}), especially by the Stieltjes inversion formula
 \bb
  \label{Psi}
 \psi (t) - \psi (0) =-\frac1{\pi} \lim_{y \to 0^+} \int_0^t \Im F(x+i y;L) dx,
 \ee
we conclude that holds
\begin{equation}
 \label{F}
\mathcal{F}(z;L)=\int F(z;L)dz=\frac14
\Bigl[z^2-2Lz-(z-L+1)R(z;L)-l_1(z)+l_2(z)\Bigr],
\end{equation}
where
 \begin{eqnarray*}
  l_1(z)&=&2(3L+1)\log\Bigl[z-(L+1)+R(z;L)\Bigr] \\
 l_2(z)&=&2(L-1)\log\Bigl[\frac{-(L-1)R(z;L)-(L-1)^2+z(L+1)}{z^2(L-1)^3}\Bigr]
  \end{eqnarray*}
Rewriting the function $R(z;L)$ in the form
$$
R(z;L)=\sqrt{(z-L-1)^2-4L}
$$
and replacing $z=x+iy$, we have
$$
R(x;L)=\lim_{y \to 0^+}R(x+iy;L)=\left\{
\begin{array}{ll}
    i\sqrt{4L-(x-L-1)^2}, \quad x \in (a,b);\\
    \sqrt{(x-L-1)^2-4L},  \quad otherwise,  \\
    \end{array}%
\right.
$$
where
\bb
a=(\sqrt{L}-1)^2,\qquad b=(\sqrt{L}+1)^2.
\ee

In the case when $x \notin \Bigl((\sqrt{L}-1)^2,(\sqrt{L}+1)^2 \Bigr)$, value $R(x;L)$ is real.
Therefore we can calculate imaginary part of $\mathcal{F}(x;L)=\lim_{y \to
0^+}\mathcal{F}(x+iy;L)$:
\begin{eqnarray*}
\Im \mathcal{F}(x;L)=\Im [l_2(x)-l_1(x)]=0.
\end{eqnarray*}

Otherwise, if $x \in \Bigl((\sqrt{L}-1)^2,(\sqrt{L}+1)^2 \Bigr)$ we have that:
\begin{eqnarray*}
l_1(x)&=2(3L+1)\log\Bigl[x-(L+1)\pm i\sqrt{4L-(x-L-1)^2}\Bigr]\\
\Im l_1(x)&=\left\{%
\begin{array}{ll}
    2(3L+1) \arctan \frac{\sqrt{4L-(x-L-1)^2}}{x-(L+1)}, & x \geq L+1; \\
    2(3L+1) \Bigl( \pi+ \arctan \frac{\sqrt{4L-(x-L-1)^2}}{x-(L+1)} \Bigr), & x<L+1 \\
\end{array}%
\right.\\
l_2(x)&=2(L-1)\log\Bigl[\frac{-(L-1)^2+2x(L+1)-i(L-1)\sqrt{4L-(x-L-1)^2}}{x^2(L-1)^3}\Bigr]\\
\Im l_1(x)&=\left\{%
\begin{array}{ll}
    2(L-1) \Bigl(2\pi+ \arctan \frac{x(L+1)-(L-1)^2}{\sqrt{4L-(x-L-1)^2}} \Bigr), & x \geq \frac{(L-1)^2}{L+1}; \\
    2(L-1) \Bigl( \pi+ \arctan \frac{x(L+1)-(L-1)^2}{\sqrt{4L-(x-L-1)^2}} \Bigr), & x<\frac{(L-1)^2}{L+1} \\
\end{array}%
\right.
\end{eqnarray*}
After substituting all considered cases in (\ref{F}), we finally
obtain the value
$$\Im \mathcal{F}(x;L)=\lim_{y \to 0^+} \Im \mathcal{F}(x+iy;L)=\Im l_2(x)-\Im l_1(x)-(x-L+1)\sqrt{4L-(x-L-1)^2}$$
From the relation ( \ref{Psi}), we conclude that
\begin{equation}
\omega(x;L)=\psi'(x;L)=-\frac1{\pi} \frac{d}{dx} \Im \mathcal{F}(x;L)
\end{equation}
and finally, we obtain
\bb   \label{wei}
\omega(x;L)=\frac1{2\pi}(1+\frac1x)\sqrt{4L-(x-L-1)^2}=\frac{\sqrt{L}}{\pi}\Bigl(1+\frac1x
\Bigr) \sqrt{1-\Bigl(\frac{x-L-1}{2\sqrt{L}}\Bigr)^2}
\ee
Previous formula holds for $x\in (a, b)$, and otherwise is $\omega(x;L)=0$.

%\section 4
\section{Determining the three-term recurrence relation}

The crucial moment in our proof of the conjecture is to determine the sequence of polynomials $\{Q_n(x)\}$
orthogonal with respect to the weight $w(x;L)$ given by
(\ref{wei}) on the interval $(a, b)$ and to find the sequences
$\{\alpha_n\}$ $\{\beta_n\}$ in the three-term recurrence relation.

%{\bf Example 4.1}

 \begin{exm}\rm
For $L=4$, we can find the first members
$$
\aligned
Q_0(x) &= 1,\qquad \qquad \qquad   \qquad\qquad   \qquad \ \|Q_0\|^2 = 5,\\
Q_1(x) &= x-\frac{24}5, \qquad \qquad  \qquad \qquad   \qquad\ \|Q_1\|^2 = \frac{104}{5},\\
Q_2(x) &=  x^2 -\frac{127}{13}x+\frac{256}{13},\qquad \qquad   \qquad\|Q_2\|^2 =\frac{1088}{13},\\
Q_3(x)&=x^3-\frac{541}{17}x^2+\frac{1096}{17}x-\frac{1344}{17},\quad \|Q_3\|^2 = \frac{5696}{17},
\endaligned
$$
wherefrom
$$
\alpha_0=\frac{24}5, \quad \beta_0=5, \qquad
\alpha_1=\frac{323}{65}, \quad \beta_1=\frac{104}{25}, \qquad
\alpha_2=\frac{1104}{221}, \quad \beta_2=\frac{680}{169}.
$$
Hence
$$
h_1=a_0=5, \quad h_2=a_0^2\beta_1\ =104, \qquad
h_3=a_0^3\beta_1^2\beta_2=5^3\ \Bigl(\frac{104}{25}\Bigr)^2
\frac{680}{169}=8704.
$$
\end{exm}

At the beginning, we will notice that in the definition of the
weight function appears the square root member.

That's why, let us consider the monic orthogonal polynomials
$\{S_n(x)\}$ with respect to the $p^{(1/2,1/2)}(x)=\sqrt{1-x^2}$
on the interval $(-1,1)$. These polynomials are monic Chebyshev
polynomials of the second kind:
$$
S_n(x)=\frac{\sin \bigl((n+1)\arccos x\bigr)}{2^n\cdot
\sqrt{1-x^2}}
$$
They satisfy the three-term recurrence relation (Chihara
\cite{chihara}):
\bb
 \label{Cheb}
S_{n+1}(x)=(x -\alpha^*_n) \ S_n(x) - \beta^*_n S_{n-1}(x) \quad
(n=0,1,\ldots), \ee with initial values
$$
S_{-1}(x)=0, \qquad S_0(x)=1,
$$
where
$$
\alpha^*_n=0\quad (n \ge 0) \qquad \text{and} \qquad
\beta^*_0=\frac{\pi}{2}, \quad \beta^*_n=\frac14 \qquad (n \ge 1).
$$
If we use the weight function $\hat w(x)=(x-c)\ p^{(1/2,1/2)}(x),$
then the corresponding coefficients $\hat \alpha_n$ and $\hat
\beta_n$ can be evaluated as follows (see, for  example,  Gautschi
\cite{gautschi})
\bb
 \label{Rec}
\aligned
\lambda_n&=S_n(c), \\
\hat \alpha_n &= c- \frac{\lambda_{n+1}}{\lambda_n}-\beta^*_{n+1} \frac{\lambda_{n}}{\lambda_{n+1}},\\
\hat \beta_n  &=\beta^*_n
\frac{\lambda_{n-1}\lambda_{n+1}}{\lambda^2_n} \qquad \qquad \quad
(n\in \mathbb N_0).
\endaligned
\ee
From the relation (\ref{Cheb}), we conclude that the sequence
$\{\lambda_n\}_{n \in \mathbb{N}}$ satisfies the following
recurrence relation:
\begin{equation}
 \label{E}
4\lambda_{n+1}-4c\lambda_n + \lambda_{n-1}=0 \qquad
(\lambda_{-1}=0; \ \lambda_0=1).
\end{equation}
The characteristic equation
$$
4z^2-4\; c\ z+1=0
$$
has the solutions
$$
z_{1,2} = \frac12\Bigl(c\pm \sqrt{c^2-1}\Bigr).
$$
and the integral solution of (\ref{E}) is
$$
\lambda_n = E_1 z_1^n + E_2 z_2^n \qquad  (n\in\mathbb N).
$$
We evaluate values $E_1$ and $E_2$ from the initial conditions
$(\lambda_{-1}=0; \ \lambda_0=1)$.

 In other to solve our problem, we will
choose $c=-\frac{L+2}{2\sqrt{L}}$. Hence
$$
z_k=\frac{-t_k}{4\sqrt{L}}\quad (k=1,2), \quad {\rm where} \quad
t_{1,2}=L+2\pm\sqrt{L^2+4}.
$$
Finally, we obtain:
$$
\lambda_n=\frac{(-1)^n}{2\cdot 4^n
L^{\frac{n}2}\sqrt{L^2+4}}\Bigl(t_1^{n+1}-t_2^{n+1}\Bigr)\quad
(\lambda=-1,0,1,\ldots),
$$
i.e,
$$
\lambda_n=\frac{(-1)^n}{2\cdot 4^n L^{\frac{n}2}\xi}\
\psi_{n+1}\quad (\lambda=-1,0,1,\ldots).
$$
 After replacing in (\ref{Rec}), we obtain:
\begin{eqnarray}
\hat \alpha_n&=&-\frac{L+2}{2\sqrt{L}} + \frac1{4\sqrt{L}}\cdot
\frac{\psi_{n+2}}{\psi_{n+1}} +
\sqrt{L}\cdot \frac{\psi_{n+1}}{\psi_{n+2}},\\
\hat \beta_n&=&\frac{\psi_{n}\psi_{n+2}}{4\psi_{n+1}^2}.
\end{eqnarray}
If  a new weight function $\tilde w(x)$ is introduced by
$$
\tilde w(x)=\hat w(ax+b)
$$
then we have
$$
\tilde \alpha_n= \frac{\hat\alpha_n-b}{a}, \qquad
\tilde \beta_n=\frac{\hat\beta_n}{a^2}\qquad (n\geq 0).
$$
 Now, by using $x\mapsto \frac{x-L-1}{2\sqrt{L}}$,\ i.e.,
 $a=\frac{1}{2\sqrt{L}}$ and $b=-\frac{L+1}{2\sqrt{L}}$, we have
the weight  function
$$
\tilde w(x)=\hat w(\frac{x-L-1}{2\sqrt{L}})=\frac12 \Bigl(\frac{x-L-1}{2\sqrt{L}}+\frac{L+2}{2\sqrt{L}}\Bigr)
\sqrt{1-\Bigl( \frac{x-L-1}{2\sqrt{L}} \Bigr)^2}.
$$
Thus
\begin{eqnarray}
\tilde \alpha_n=-1 + \frac1{2}\cdot \frac{\psi_{n+2}}{\psi_{n+1}}
+ 2L\cdot \frac{\psi_{n+1}}{\psi_{n+2}} \qquad (n\in\mathbb N_0),
\end{eqnarray}

and

\begin{eqnarray}
 \tilde \beta_0=(L+2)\frac{\pi}2,\qquad
\tilde\beta_n=L\frac{\psi_{n}\psi_{n+2}}{\psi_{n+1}^2} \quad
(n\in\mathbb N).
\end{eqnarray}

\begin{exm}\rm
For $L=4$, we get
$$
\aligned
P_0(x) &= 1,\qquad \qquad \qquad   \qquad\qquad   \qquad \ \|P_0\|^2 = 3\pi,\\
P_1(x) &= x-\frac{17}3, \qquad \qquad  \qquad \qquad   \qquad\ \|P_1\|^2 = \frac{32\pi}{3},\\
P_2(x) &=  x^2 -\frac{43}{4}x+\frac{101}{4},\qquad \qquad   \qquad\|P_2\|^2 =42\pi,\\
P_3(x)&=x^3-\frac{331}{21}x^2+\frac{1579}{21}x-\frac{2189}{21},\quad
\|P_3\|^2 = \frac{3520\pi}{21},
\endaligned
$$
wherefrom
$$
\tilde \alpha_0=\frac{17}3, \quad \tilde \beta_0=3\pi, \qquad
\tilde \alpha_1=\frac{61}{12}, \quad \tilde \beta_1=\frac{32}{9},
\qquad \tilde \alpha_2=\frac{421}{84}, \quad
\tilde\beta_2=\frac{63}{16}.
$$

\end{exm}

Introducing the weight
$$
\breve{w}(x)=\frac{2L}{\pi} \tilde w(x)
$$
will not change the monic polynomials and their recurrence
relations, only it will multiply the norms by the factor $2L/\pi$,
i.e.
$$
\breve{P}_k(x)\equiv P_k(x), \qquad
 \|\breve{P}_k\|^2_{\breve{w}}=\int_a^b
 \breve{P}_k(x)\breve{w}(x)\ dx = \frac{2L}{\pi} \|P_k\|^2 \qquad
 (k\in\mathbb N_0),
$$

$$
\breve{\beta_0} =L(L+2),
 \quad \breve{\beta_k} = \tilde \beta_k \quad (k\in\mathbb N),
\qquad \qquad\qquad  \breve{\alpha_k} = \tilde\alpha_k
\qquad(k\in\mathbb N_0).
$$
Here is
\begin{eqnarray}
\breve{\beta_0}\breve{\beta_1}\cdots {\breve
\beta}_{n-1}=\frac{L^{n}}2\cdot\frac{\psi_{n+1}}{\psi_{n}}.
\end{eqnarray}

In the book \cite{wgautschi}, W. Gautschi has treated the next
problem: If we know all about the MOPS orthogonal with respect to
$\breve{w}(x)$ what can we say about the sequence $\{Q_n(x)\}$
orthogonal with respect to a weight
$$
w_d(x)=\frac{\breve{w}(x)}{x-d}\qquad (d\notin
\text{support}(\tilde w))\ ?
$$
W. Gautshi has proved that, by the auxiliary sequence
$$
r_{-1} = -\int_{\mathbb R} w_d(x)\ dx, \qquad r_n = d-\breve{
\alpha}_n-\frac{\breve{\beta}_n}{r_{n-1}} \quad (n=0,1,\ldots),
$$
it can be determined
$$
\aligned \alpha_{d,0}&=\breve{\alpha}_0+r_0,\qquad
&\alpha_{d,k}&=&\breve{\alpha}_k+r_k-r_{k-1}&,\\
\beta_{d,0}&=-r_{-1}, \qquad &\beta_{d,k}&=&\breve{
\beta}_{k-1}\frac{r_{k-1}}{r_{k-2}}\qquad & \qquad (k\in\mathbb
N).
\endaligned
$$
In our case it is enough to take $d=0$ to get the final weight
$$
w(x)= \frac{\breve{w}(x)}{x}.
$$
Hence
\begin{eqnarray}
r_{-1} = -(L+1), \qquad r_n =
-\Bigl(\breve{\alpha}_n+\frac{\breve{ \beta}_{n}}{r_{n-1}}\Bigr)
\quad (n=0,1,\ldots).
\end{eqnarray}

\begin{lem} The parameters $r_n$ have the explicit form
\begin{eqnarray}
r_n= - \frac{\psi_{n+1}}{\psi_{n+2}}\cdot
\frac{L\psi_{n+2}+\xi\varphi_{n+2}}{L\psi_{n+1}+\xi\varphi_{n+1}}\qquad
(n\in\mathbb N_0).
\end{eqnarray}
\end{lem}
\noindent{\it Proof.} We will use the mathematical induction. For
$n=0$, we really get the expected value
$$
r_0=-\frac{L^2+2L+2}{(L+1)(L+2)}.
$$
Suppose that it is true for $k=n$. Now,  by the properties for $\varphi_n$ and $\psi_n$, we have
$$
\tilde\alpha_{n+1}\cdot r_n +
\tilde\beta_{n+1}
=-\frac{\psi_{n+1}}{\psi_{n+3}}\cdot
\frac{L\psi_{n+3}+\xi\varphi_{n+3}}{L\psi_{n+1}+\xi\varphi_{n+1}}.
$$
Dividing with $r_n$, we conclude that the formula is valid for $r_{n+1}$. $\Box$

\smallskip

\begin{exm}\rm
For $L=4$, we get
$$
r_{-1}=-5,\quad r_{0}=-\frac{13}{15},\quad
r_{1}=-\frac{51}{52},\quad r_{2}=-\frac{356}{357},
$$
wherefrom
$$
 \alpha_0=\frac{24}5, \quad \beta_0=5, \qquad
\alpha_1=\frac{323}{65}, \quad \beta_1=\frac{104}{25}, \qquad
\alpha_2=\frac{1104}{221}, \quad \beta_2=\frac{680}{169},
$$
just the same as in the Example 4.1.
\end{exm}

{\it Proof of the main result}. The Krattenthaler's formula
(\ref{formula}) can be also written in the form
  \bb
 \label{formula1}
 h_1=a_0,\qquad h_n= \beta_{0}\beta_{1} \beta_{2}\cdots \beta_{n-2} \beta_{n-1} \cdot h_{n-1}.
 \ee

 From the theory of orthogonal polynomials, it is known that
 \bb
 \label{formula2}
 \|Q_{n-1} \|^2 = \beta_{0}\beta_{1} \beta_{2}\cdots \beta_{n-2} \beta_{n-1}\quad
 (n=2,3,\ldots),
 \ee
wherefrom
  \bb
 \label{formula3}
  h_1=a_0,\qquad h_n= \|Q_{n-1} \|^2 \cdot h_{n-1}\quad (n=2,3,\ldots).
 \ee
 Here,
  \bb
 \label{formula4}
 \|Q_{n-1} \|^2 = \beta_{0}\frac{r_{n-2}}{r_{-1}}\prod_{k=0}^{n-2}\breve\beta_{k}
 =\frac{L^{n-1}}2\cdot\frac{L\psi_{n}+\xi\varphi_{n}}{L\psi_{n-1}+\xi\varphi_{n-1}} .
 \ee

We will apply the mathematical induction again. The formula for
$h_n$ is true for $n=1$. Suppose that it is valid for $k=n-1$.
Then
$$
h_n=
\frac{L^{n-1}}2\cdot\frac{L\psi_{n}+\xi\varphi_{n}}{L\psi_{n-1}+\xi\varphi_{n-1}}
\cdot \frac{L^{(n-1)(n-2)/2}}{2^{n}\xi}\cdot \left(L\psi_{n-1}+\xi
\varphi_{n-1}\right),
$$
wherefrom  it follows that the final statement
$$
h_n= \frac{L^{n(n-1)/2}}{2^{n+1}\xi}\cdot \left(L\psi_n+\xi
\varphi_n\right)\qquad (n\in \mathbb N)
$$
is true. $\Box$

\smallskip

{\bf  Predrag Rajkovi\'c, Marko D. Petkovi\'c,}\\
 University of Ni\v s, Serbia and Montenegro  \\
Address: A. Medvedeva 14, 18000 Ni\v s, Serbia and Montenegro\\
e-mail: {\tt pedja.rajk@gmail.com, dexterofnis\@gmail.com }

\smallskip

{\bf Paul Barry} \\
School of Science, Waterford Institute of Technology, Ireland\\
e-mail: {\tt pbarry@wit.ie}\\

\end{document}